\documentclass[11pt,leqno]{amsart}
\usepackage{amsthm,amsfonts,amssymb,amsmath,oldgerm}
\numberwithin{equation}{section}
\usepackage{fullpage}
\usepackage{amsmath}
\usepackage{amsfonts}
\usepackage{amssymb}
\usepackage{setspace}
\usepackage{stmaryrd}
\usepackage{mathrsfs}
\usepackage{fancyhdr}
\usepackage{graphicx}

\usepackage{psfrag}
\usepackage{listings} 

\usepackage[top=30mm,bottom=30mm,left=25mm,right=25mm,a4paper]{geometry}

\renewcommand\d{\partial}
\newcommand\dD{\textrm{d}}
\newcommand\eD{\textrm{e}}
\newcommand\iD{\textrm{i}}
\newcommand\ii{\iD}

\newcommand{\Id}{{\rm Id}}

\newcommand\br{\begin{remark}}
\newcommand\er{\end{remark}}
\newcommand\bp{\begin{pmatrix}}
\newcommand\ep{\end{pmatrix}}
\newcommand{\be}{\begin{equation}}
\newcommand{\ee}{\end{equation}}
\newcommand\ba{\begin{equation}\begin{aligned}}
\newcommand\ea{\end{aligned}\end{equation}}

\newcommand{\beg}{\begin{example}}
\newcommand{\eeg}{\end{exaplem}}
\newcommand{\bpr}{\begin{proposition}}
\newcommand{\epr}{\end{proposition}}
\newcommand{\bt}{\begin{theorem}}
\newcommand{\et}{\end{theorem}}
\newcommand{\bc}{\begin{corollary}}
\newcommand{\ec}{\end{corollary}}
\newcommand{\bl}{\begin{lemma}}
\newcommand{\el}{\end{lemma}}
\newcommand{\bd}{\begin{definition}}
\newcommand{\ed}{\end{definition}}
\newcommand{\brs}{\begin{remarks}}
\newcommand{\ers}{\end{remarks}}

\newtheorem{theorem}{Theorem}[section]
\newtheorem{proposition}[theorem]{Proposition}
\newtheorem{corollary}[theorem]{Corollary}
\newtheorem{lemma}[theorem]{Lemma}

\theoremstyle{remark}
\newtheorem{remark}[theorem]{Remark}
\theoremstyle{definition}
\newtheorem{definition}[theorem]{Definition}

\newtheorem{example}[theorem]{Example}

\newcommand\R{\mathbf R}
\newcommand\C{\mathbf C}
\newcommand{\N}{\mathbf N}
\newcommand{\Z}{\mathbf Z}

\newcommand{\intg}{[\![}
\newcommand{\intd}{]\!]}

\newcommand\bD{{\mathbf D}}

\newcommand\bF{{\mathbf F}}

\newcommand\bM{{\mathbf M}}

\newcommand\bU{{\mathbf U}}
\newcommand\bV{{\mathbf V}}

\newcommand\bff{{\mathbf f}}

\newcommand\ubM{{\underline \bM}}

\newcommand\ubU{{\underline \bU}}

\newcommand\uU{{\underline U}}

\newcommand\uk{{\underline k}}

\newcommand\uom{{\underline \omega}}

\newcommand\cM{{\mathcal M}}
\newcommand\cN{{\mathcal N}}

\newcommand\tW{\widetilde{W}}

\title{
Space-modulated stability and averaged dynamics
}

\author{L.Miguel Rodrigues}
\address{
Universit\'e de Rennes 1,
IRMAR, UMR CNRS 6625,
263 avenue du General Leclerc;
F-35042 Rennes Cedex, FRANCE}
\email{{\tt luis-miguel.rodrigues@univ-rennes1.fr}}
\thanks{Research of L.Miguel Rodrigues was partially supported by the ANR project
BoND ANR-13-BS01-0009-01.\\
}

\begin{document}

\begin{abstract}
In this brief note we give a brief overview of the comprehensive theory, recently obtained by the author jointly with Johnson, Noble and Zumbrun, that describes the nonlinear dynamics about spectrally stable periodic waves of parabolic systems and announce parallel  results for the linearized dynamics near cnoidal waves of the Korteweg--de Vries equation. The latter are expected to contribute to the development of a dispersive theory, still to come.
\end{abstract}

\date{\today}
\maketitle

%\begin{center}
{\it Keywords}: periodic traveling waves; stability; modulation.
%\end{center}

%\begin{center}
{\it 2010 MSC}: 35B10, 35B35, 35K59, 35P05, 35Q53, 37K35.%?

%%%%%%%%%%%%%%%%%%%%%%%%%%%%%%%%%%%%%%%%%%%%%%%%%%%%%% 
%\end{center}

%\clearpage
\tableofcontents
%\clearpage

%%%%%%%%%%%%%%%%%%%%%%%%%%%%%%%%%%%%%%%%%%%%%%%%%%%%%%         
%       INTRODUCTION                                 %
%%%%%%%%%%%%%%%%%%%%%%%%%%%%%%%%%%%%%%%%%%%%%%%%%%%%%%

\section{Introduction}\label{s:introduction}

A classical approach to the qualitative analysis of dynamical behavior driven by partial differential equations starts with the identification of special classes of solutions that exhibit some simple internal structure. One hope is that when these coherent structures possess a nontrivial basin of attraction they could serve as elementary blocks to describe a richer complex dynamics.

However, in many situations, determining dynamical stability even for the simplest reference solutions appears as a daunting task. Hence the general strategy consisting in obtaining on one hand general theorems --- or at least general accurate frameworks --- to derive dynamical behavior from spectral information on linearized evolution and on the other hand in gathering separately the required pieces of spectral information. Fortunately it turns out that in many cases the latter may at least be obtained either analytically in some asymptotic regimes or numerically for (a large part of) the full family of background solutions. Interestingly enough, besides those classical ways of determining spectral stability, an approach based on interval arithmetics and leading to computer assisted proofs of spectral stability seems to have reached sufficient maturity to completely solve some stability issues; the reader is referred to \cite{Barker} for an instance of such a line of study. 

With this general scheme in mind some considerable amount of effort has been devoted recently and is still directed towards the development of a general stability theory for periodic traveling wave solutions of systems of partial differential equations. The purpose of this note is to offer both a terse introduction to the field and an overview of the recent activity of the author, partly jointly with others. 

As a preliminary warning we want to strongly emphasize that traveling waves under consideration are solutions of equations modeling extended systems. In particular our physical domains are unbounded, at least in the direction of propagation. This situation arises naturally when idealizing a phenomena where boundaries seems to play a minor role as in the coherent motion of waves over large distances. As a consequence in our set of problems periodicity in space is a feature of the internal structure of background waves and not of the physical domain.

From the point of view of stability this results in the fact that we are not interested in co-periodic or sub-harmonic perturbations but in stability under arbitrary smooth and localized ones. Consequently we aim also at taking into account the rich multiscale spatio-temporal dynamics emerging from the slow modulations of periodic waves. In contrast with the local analysis around patterns with a simpler internal spatial structure  such as solitary waves, fronts, kinks or shocks, and even for simpler-looking systems as those modeling reaction-diffusion processes, this inherent complexity has precluded until very recently any completely satisfactory mathematical analysis of dynamical stability and large-time asymptotic behavior. 

Nevertheless, for general parabolic systems, after the pioneering work of Schneider \cite{Schneider-SH,Schneider-proc} at the end of the past century, an intense activity in the last few years --- more or less culminating in \cite{JNRZ-conservation} --- has led to a clear picture of the relevant notion of spectral stability, called diffusive spectral stability, and of the asymptotic description of the dynamics near stable periodic waves. The general theory even includes systems that are only parabolic in the pointwise hypocoercive sense of Kawashima, thus effectively including at least small waves of viscous compressible models, and have been extended recently to some systems that are only symmetrizable in a suitable averaged sense \cite{RZ}, in particular relaxing the smallness condition for roll-waves of the St. Venant system describing surface waves on fluid film flowing down an incline. Incidentally, we point out that the diffusive spectral stability of the latter roll-waves have been extensively studied by a combination of numerical methods and spectral perturbation analysis \cite{BJNRZ-KdV-SV,BJNRZ-KdV-SV-note}, leading --- jointly with the abstract nonlinear theory --- to a fairly comprehensive description of the dynamics near those primary instabilities.

In contrast, for dispersive equations, as far as the author knows, even a single example of a proof of nonlinear stability is missing. Yet the detailed description of the \emph{linear} dynamics of solutions to the Korteweg--de Vries near cnoidal waves --- recently obtained in \cite{R_linKdV} --- shades some light on what kind of nonlinear dynamics is to be expected.

Based essentially on \cite{JNRZ-conservation,R_linKdV} our goal here is to give both an overview of the well-developed parabolic theory and a glimpse at what could be expected of a dispersive theory, still to come. On the former the reader may also benefit from consulting \cite{R}. We warn the reader that, to keep the presentation as short as possible, in most of the present note the general tone is rather descriptive and informal.

\section{General information}

We shall restrict our discussion to planar traveling waves --- as these are the objects for which the theory is more advanced --- and we shall regard them as one-dimensional objects. The latter choice comes with two main consequences. On one hand nonlinear analysis is considerably harder in space dimension one since scattering mechanisms are weaker thus decay resulting from diffusive-like or dispersive-like behavior is slower\footnote{Think about decay rate $t^{-\tfrac{N}{2}}$ in $L^\infty$-norms for solutions to the heat or the Schr\"odinger equation in space-dimension $N$ starting from integrable initial data.}. This is reflected by the fact that the first results of nonlinear stability of spectrally stable periodic waves of parabolic systems of conservation laws were first proved in dimension higher than two, then in dimension two and finally in dimension one \cite{Oh-Zumbrun-cons-low-freq-multiD,Oh-Zumbrun-cons-nonlinear,Oh-Zumbrun-cons-nonlinear-erratum,Johnson-Zumbrun-cons-generic}, \cite{Johnson-Zumbrun-cons-1D-2D}, \cite{JNRZ-conservation}. This is also responsible for the asymptotically \emph{linear} behavior exhibited by localized perturbations of periodic waves in dimension higher than one. On the other hand, spectral stability may seem easier to met in dimension one since in some sense spectral stability in dimension one may be interpreted as a higher-dimensional spectral stability but restricted to co-planar perturbations. Obviously one may then rightfully wonder what could be the physical meaning of a situation where a planar wave would be stable as a one-dimensional object but unstable as a higher-dimensional object. An element of answer lies in the fact that our one-dimensional solutions sometimes correspond in higher dimension not to planar waves but to genuinely multidimensional objects resulting from some confinement mechanisms that may preserve stability\footnote{Think about a fluid flowing down a channel.}. In short, dimension one offers both a richer nonlinear dynamics --- at the cost of a substantially more involved analysis --- and better hope to find stable waves.

Moreover, in this brief note we discuss only two sets of equations, on one side the Korteweg--de Vries equation (KdV)
\be\label{KdV}
U_t\ +\ \left(\tfrac12\,U^2\right)_x\ +\ U_{xxx}\ =\ 0\,.
\ee
and, on the other side, a general second-order semi-linear parabolic system of conservation laws
\be\label{conservation}
\bU_t\,+\,(\bff(\bU))_x\,=\,\bD\,\bU_{xx}\,.
\ee
In the above, time variable $t$ belongs either to $\R$ or $\R_+$ and spatial variable $x$ belongs to $\R$. The unknown of \eqref{KdV} is scalar, $U(t,x)\in\R$, whereas the unknown of \eqref{conservation} is vector-valued, $\bU(t,x)\in\R^d$, for some $d\in\N^*$, the corresponding advection flux is given by a smooth function $\bff:\R^d\to\R^d$ and the diffusion matrix $\bD$ is some $d\times d$ positive definite symmetric matrix.

To start being more specific let us add that here by a periodic wave we mean a special solution $U(t,x)=\uU(kx+\omega t)$ given by a periodic profile $\uU$ and traveling with uniform velocity $c=-\omega/k$. When normalizing the period of $\uU$ to one --- as we shall do --- $k$ accounts for spatial wavenumber and $\omega$ for time frequency.

The choice of \eqref{conservation} is purely expository as it does not bring any significant simplification in statements and proofs, which apply almost word-by-word to a much larger class of parabolic systems. In contrast, the choice of \eqref{KdV} is commanded by very specific reasons. Indeed, prior to \cite{R_linKdV}, the author was ignorant of what could be the best notion of spectral "dispersive" stability that one may expect and which kind of spectral conditions could ensure some form of linear asymptotic stability. Hence the choice of \eqref{KdV} that, beyond its own importance as a universal amplitude equation for long weakly-nonlinear phenomena, enjoys the advantage of being integrable so that in principle everything on its periodic waves and their spectra may be explicitly derived. We stress however that the main intention of \cite{R_linKdV} is not to study the Korteweg--de Vries equation specifically but to develop a strategy, first illustrated on \eqref{KdV}, that transfers spectral knowledge to linear bounds. This opens possible applications to cases where the spectrum is not known explicitly but where required spectral conditions may be either proved in an asymptotic regime, for instance for small waves, or investigated numerically.

Since we aim at discussing dynamics near a specific wave described by some profiles $\uU$ and with wavenumber $\uk$ and frequency $\uom$ it is convenient\footnote{But not absolutely mandatory at least for some of the arguments. See for instance \cite{KR} where translational invariance is essentially broken by discreteness of the spatial domain and no obvious way to turn traveling waves in steady solutions is available.} to switch to the moving frame of this specific wave. For instance, defining implicitly $W$ through $U(t,x)\ =\ W(t,\uk\,x+\uom\,t)$, equation \eqref{KdV} becomes
\be\label{KdV-move}
W_t\ +\uom\,W_x\ +\ \uk\,\left(\tfrac12\,W^2\right)_x\ +\ \uk^3\,W_{xxx}\ =\ 0
\ee
and $\uU$ is a stationary solution of \eqref{KdV-move}. A (naive) linearization of \eqref{KdV-move} about $\uU$ leads to the consideration of $\tW_t-L\tW=0$ where
\be\label{linop}
L\,\tW\ =\ -\,\uom\,\tW_x\ -\ \uk\,\left(\uU\,\tW\right)_x\ -\uk^3\,\tW_{xxx}\,.
\ee

\subsection{The Bloch transform}

This naturally suggests the introduction of an integral transform specially tailored to analyze differential operators with periodic coefficients, here normalized to period one.  To motivate its definition, let us observe that in this context the simplest relevant class of functions that may span the full set of arbitrary functions is given by so-called Bloch waves that have a simple behavior under translation by one period as they are simply multiplied by some number of modulus one, say $\eD^{\iD\xi}$ where $\xi$ is a Floquet exponent. Equivalently those are the functions that may be written as $x\mapsto\eD^{\iD\xi\,x}\check{u}(x)$ for some Floquet exponent $\xi\in[-\pi,\pi)$ and some function of period one $\check{u}$.

One rationale behind the definition of the \emph{Bloch transform} --- sometimes called Floquet-Bloch transform --- is that one may indeed achieve a decomposition 
\be\label{inverse-Bloch}
g(x)\ =\ \int_{-\pi}^\pi \eD^{\ii\xi x}\ \check{g}(\xi,x)\ \dD\xi,
\ee
with each $\check{g}(\xi,\cdot)$ periodic of period one, by summing modes of the Fourier decomposition sharing the same Floquet exponent
\be\label{Bloch}
\check{g}(\xi,x)\ :=\ \sum_{j\in\Z} \eD^{\ii\,2j\pi x}\ \widehat{g}(\xi+2j\pi).
\ee
where
$$
\widehat g(\xi)\ :=\ \frac{1}{2\pi} \int_\R \eD^{-\ii\xi x} g(x)\ \dD x.
$$
An equivalent approach consists in introducing a periodization of $x\mapsto\eD^{-\ii\xi x}g(x)$
$$
\check{g}(\xi,x)\ :=\ \sum_{\ell\in\Z}\eD^{-\ii\xi (x+\ell)}g(x+\ell)
$$
which is readily seen to coincide with the above formula through the Poisson summation formula.

Bloch analysis has essentially the same flavor as classical Fourier analysis\footnote{That is why we have been rather loose with convergence issues in the above definitions since they are solved in the classical way.}. Observe for instance that --- modulo a multiplicative constant --- the Bloch transform is a total isometry between $L^2(\R)$ and $L^2((-\pi,\pi),L_{per}^2((0,1)))$ and that interpolating between Parseval identity and triangle inequalities leads to the following Hausdorff-Young inequalities, for $2\leq p\leq\infty$,
$$
\|g\|_{L^p(\R)}\ \leq\ (2\pi)^{1/p}
\| \check{g}\|_{L^{p'}([-\pi,\pi],L^p((0,1)))}\,,\qquad
\|\check{g}\|_{L^{p}([-\pi,\pi],L^{p'}((0,1)))}\ \leq
\ (2\pi)^{-1/p}\|g\|_{L^{p'}(\R)}\,
$$
where $p'$ denotes conjugate Lebesgue exponent, $1/p+1/p'=1$. Going on with the former point it may be rather instructive\footnote{This is probably the point of view that generalizes in the most straight-forward way to the spatially discrete case, where, for any $N\in \N$, a suitable discrete Bloch transform provides an isometry between $\ell^2(\Z)$ and $L^2((-\pi,\pi),\ell_{per}^2(\intg0,N\intd))$. See \cite{KR}.} to regard $g\in L^2(\R)$ as an element of $\ell^2(\Z,L^2((0,1)))$ providing the $L^2((0,1))$-valued Fourier series of $\xi\mapsto \eD^{\ii\xi \cdot}\check{g}(\xi,\cdot)\in L^2((-\pi,\pi),L^2((0,1)))$.

As expected, differential operators $L$ with coefficients of period one and acting on a space of functions defined on the whole line, say $L^2(\R)$, may then be thought on the Bloch side as operator-valued multipliers acting Floquet exponent by Floquet exponent, that is
$(Lg)\,\check{ }\,(\xi,\cdot)=L_\xi \check{g}(\xi,\cdot)$ where for each $\xi\in[-\pi,\pi)$ the operator $L_\xi$ acts on a space of functions of period one, say $L_{per}^2((0,1))$. This follows readily from applying $L$ to \eqref{inverse-Bloch}, a process that also gives 
$$
L_\xi\ =\ \eD^{-\ii\xi \cdot}\ L\ \eD^{\ii\xi \cdot}
$$
where $L$, $L_\xi$ are here used as formal operators.

This "diagonalization" has direct counterparts for spectra and generated semi-groups. Observe in particular that for operators $L$ obtained from linearization (in a suitable frame) of either \eqref{KdV} or \eqref{conservation} each $L_\xi$ has compact resolvents (when acting say on $L_{per}^2((0,1))$) hence has a spectrum made of a discrete set of eigenvalues with finite multiplicity and the associated spectral map $\xi\mapsto \sigma(L_\xi)$ is continuous. Therefore the spectrum of $L$  consists locally in curves parametrized by Floquet exponents and 
$$
\sigma_{L^2(\R)}(L)\ =\ \bigcup_{\xi\in[-\pi,\pi)}\sigma_{L_{per}^2((0,1))}(L_\xi)\,.
$$
As an illustration we mention that accordingly the minimal notion of spectral stability, that is $\sigma(L)\subset\{\lambda\,|\,\textrm{Re}(\lambda)\leq0\}$, is also written as $\sigma(L_\xi)\subset\{\lambda\,|\,\textrm{Re}(\lambda)\leq0\}$ for any $\xi\in[-\pi,\pi)$. Yet obviously such a requirement is too weak to yield any form of dynamical stability even at the linear level.

\subsection{Stability}

On the other hand, concerning the stability of traveling waves, the strongest notion of spectral stability, that requires a spectral gap $\sigma(L)\subset\{\lambda\,|\,\textrm{Re}(\lambda)\leq-\eta\}$ for some $\eta>0$, is precluded by the presence of a full family of similar traveling waves near any given wave of a specific type. Indeed at the very least translation invariance of the equations yields that any translated copy of a traveling wave is a traveling wave. Then for traveling waves with localized variations --- fronts, solitary waves... --- taking derivatives along this family of waves provides elements of the generalized kernel of $L$. However, for localized waves of parabolic systems one may still expect that $\sigma(L)\subset\{\lambda\,|\,\textrm{Re}(\lambda)\leq-\eta\}\cup\{0\}$ for some $\eta>0$ and that $0$ is an eigenvalue of finite multiplicity corresponding to the dimension of the family of neighboring traveling waves. In this case one may prove by adapting the arguments providing the central manifold theory that such waves are (orbitally) asymptotically stable and that at leading order nonlinear dynamics is reduced to a finite-dimensional evolution describing the time evolution of parameters encoding traveling waves. In short, the main dynamics occurs through modulation in time of wave parameters. Let us observe already that some changes in parameters of the waves --- for instance phase parameters, encoding the position of the wave --- have a dramatic impact on norm comparisons and the notion of dynamical stability has to be accommodated accordingly. This leads to the classical notion of orbital stability that asks for a control of the proximity of orbits of functions under the action of symmetries of the equation, which when the only relevant symmetry is space translation amounts to a bound on
$$
\inf_{\substack{\Psi \textrm{ uniform}\\\textrm{translation}}}\ \|U\circ\Psi-\uU\|_{X} 
$$
(for some chosen functional norm $\|\cdot\|_{X}$). See \cite{Henry,KapitulaPromislow_book} for some related results, references and discussions.

The situation differs significantly when periodic waves are considered. First when differentiating a family of periodic waves along a parameter independent of the period one obtains not elements of the generalized kernel of $L$ --- which is reduced to $\{0\}$ --- but of $L_0$, since obtained elements do not belong to $L^2(\R)$ but to $L^2_{per}((0,1))$. This seems to be a minor difference but then by varying $\xi$ this implies that no matter how small $\eta>0$ is chosen $\sigma(L)\setminus \{\lambda\,|\,\textrm{Re}(\lambda)\leq-\eta\}$ contains curves of spectrum passing through zero. Therefore the best one can reasonably expect is that when $\eta$ is sufficiently small then the former set is indeed reduced to those curves, that the number of such critical curves --- which is the algebraic dimension of zero as an eigenvalue of $L_0$ --- is indeed the dimension of the family of periodic waves minus one\footnote{We remind the reader that derivatives corresponding to variations of the period do not generate elements of the generalized kernel of $L_0$.} and that those curves are non degenerate at $\xi=0$. The corresponding notion of spectral stability is now well-known as \emph{diffusive spectral stability}. We recall its precise definition in Section \ref{s:parabolic}. In such a situation there is no hope for a reduction even at the leading-order to a finite-dimensional evolution that could not account for such dispersion relations. Moreover there is also no hope for some exponential decay and the decay mechanism of the critical part of the evolution is truly infinite-dimensional as it is due to scattering of solutions by a diffusive spreading. We stress that with this respect the assumption of non degeneracy of critical curves is crucial since this is the non vanishing at $\xi=0$ of the second order derivatives of the real part of spectral curves that encodes at the spectral level the diffusive character of the spectral stability. Note in particular that no matter what the order of parabolicity of the original system is the deduced decay is of heat-like type.

The issue that even a small perturbation may modify spatial position and therefore have a dramatic effect on direct norm comparison is even stronger for periodics as it occurs cell by cell for an unbounded number of cells and thus can not be compensated by some spatially uniform transformation, nor more generally by any element of a prescribed finite-dimensional set of transformations. For the sake of clarification and unification, a specially-designed notion of nonlinear stability was formalized in \cite{JNRZ-conservation} and called there \emph{space-modulated stability}. It goes as follows. Given a functional space $X$, instead of measuring at some prescribed time the proximity of a function $U$ to a background periodic wave $\uU$ directly with $\|U-\uU\|_X$, one tries to control
\be\label{sm_distance}
\delta_X(U,\uU)\ =\ \inf_{\Psi\ \textrm{one-to-one}}\quad \|U\circ\Psi-\uU\|_X\ +\ \|\d_x(\Psi-\Id)\|_X\,.
\ee
This allows, before comparison, for a resynchronization of cells through a phase change $\Psi$ that is in turn forced to remain locally in space (and time) close to uniform translations. As the notion of orbital stability does for unimodal waves, this provides a global way to quantify preservation of shapes in a periodic context. Besides being a generalization of orbital stability, it is also closely related to the Skohokhod metric on functions with discontinuities that allows for a near-identity synchronization of jumps.

As we recall in Section~\ref{s:parabolic} abstract stability results for parabolic systems, as proved in \cite{JNRZ-conservation}, show that if a periodic wave is diffusively spectrally stable it is also nonlinear stable, in a space-modulated sense, from $L^1\cap H^s$ to $H^s$ provided $s$ is large enough. However for Hamiltonian systems as \eqref{KdV}, and this even in finite dimension, diffusive spectral stability is excluded by Hamiltonian symmetry of the spectrum which implies that in a stable situation the whole spectrum lies on the imaginary axis. For finite dimensional Hamiltonian systems no decay mechanism is available and one instead prove bounded --- and not asymptotic --- stability by a Lyapunov functional argument. This actually does use a spectral gap argument but not on the dynamical generator $L$ but on the variational Hessian of the Hamiltonian. Generalizations of this argument also apply for waves with localized variations of Hamiltonian partial differential equations. But for the same reason expounded above on $L$, a spectral gap argument on the variational Hessian is precluded by the Floquet structure of its spectrum and another argument is needed --- even at the linear level --- to analyze the periodic wave problem. See however \cite{Angulo-Pava,KapitulaPromislow_book,Benzoni-Mietka-Rodrigues} for examples of what may be obtained when one restricts to co-periodic or sub-harmonic perturbations and hence enforce discreteness of Floquet exponents.

Contrary to what happens in finite dimension, note that for localized waves of Hamiltonian partial differential equations it is actually possible to go further and prove in some cases asymptotic stability by using decay provided by scattering of dispersive type. In Section~\ref{s:KdV} we show how the spreading argument is successfully generalized in \cite{R_linKdV} to prove linear asymptotic stability of cnoidal waves of \eqref{KdV}. Observe that while in the diffusive context a similar argument is used to deal with critical spectral curves passing through zero, here the whole spectrum is critical and part of the difficulty is to apply uniform arguments to spectral curves containing an infinite number of Floquet eigenvalues. We stress that linear stability also needs to be understood in a space-modulated sense that we exhibit now.

To prove nonlinear space-modulated stability of $\uU$ for \eqref{KdV-move}, one should introduce some $(V,\psi)$ related to the expected solution $W$ through 
\be\label{space-mod-ansatz}
W\circ\Psi\ =\ \uU\,+\,V\,,\qquad \Psi\ =\ \Id-\psi\,.
\ee
Stability would then be proved if one were able to choose $(W,\psi)$ in such a way that $W$ is indeed a solution and $V$ and the derivatives of $\psi$ remain small if they are sufficiently small initially. With this in mind, we observe that in terms of $(V,\psi)$ equation \eqref{KdV-move} takes the form
$$
(V+\uU_x\psi)_t\ -\ L\,(V+\uU_x\psi)\ =\ \cN[V,\psi_t,\psi_x]
$$
where $L$ is the operator defined above, and $\cN$ is nonlinear in $(V,\psi_t,\psi_x)$ and their derivatives, and, locally, at least quadratic. One important point is that the nonlinear part does not involve $\psi$ itself but only its derivatives. This is crucial since $\psi$ itself is not expected to decay but to become locally constant. Then the linearized dynamics of interest is given by
\be\label{LinKdV}
(V+\uU_x\psi)_t\ -\ L\,(V+\uU_x\psi)\ =\ 0\,.
\ee
As a result, at the linear level, one may indeed focus on the evolution generated by $L$ as we have done up to now, but instead of aiming at controlling the evolved $W$ through $\|W\|_X$ one should aim at bounding
\be\label{sm_norm}
N_X(W)\ =\ \inf_{\substack{W=V+\uU_x\psi\\\psi(\infty)=-\psi(-\infty)}}%\\\psi \textrm{ is low frequency}}} 
\|V\|_{X}\ +\ \|\psi_x\|_{X}\,.
\ee

Observe that we have added the constraint $\psi(\infty)=-\psi(-\infty)$ to ensure that $\psi=\d_x^{-1}\psi_x$, where $\d_x^{-1}$ is defined as a principal value on the Fourier side. It is intended to ensure that a knowledge of $\psi_x$ does provide some control on $\psi$. Mark that this constraint does not restrict potential applications of linear estimates to a nonlinear analysis since the centering constraint may be achieved initially by translating $\uU$ by $\tfrac12(\psi_0(\infty)+\psi_0(-\infty))$ and is then preserved by the evolution.

\subsection{Averaged dynamics}

We now turn to the question of determining a leading-order description of the dynamics. As we have already pointed out, since studied phenomena are genuinely multi-scale in time \emph{and space}, one can not expect a reduction to a finite-dimensional time evolution as for localized patterns. As we have also suggested, in dimension one --- in contrast with what happens for higher dimensional cases --- decay rates are too slow to expect that the nonlinear dynamics of \eqref{conservation} will be asymptotically linear. So we must seek another set of nonlinear partial differential equations providing an asymptotically equivalent simpler dynamics. Note also that even when the original equation is scalar, as \eqref{KdV} is, perturbations of periodic waves may split in different parts traveling at distinct characteristic velocities --- usually denominated group velocities ---, a phenomenon more easily captured by a \emph{system} of equations, so that in general\footnote{In particular this is always the case when the system contains only conservation laws as in the present note.} no reduction in the number of equations may be expected either.

What is the expected gain then ? The gain is in separation of scales. For instance, we may enumerate typical scales for \eqref{conservation} as being of size one for oscillations of the background wave, $t$ for the hyperbolic-like behavior encoded by group velocities and $\sqrt{t}$ for the heat-like diffusive character leading to diffusive spectral stability. One may hope to capture the evolution at large scales $t$, $\sqrt{t}$ by some form of averaged dynamics that does not contain any oscillation anymore. In this homogenized description the slow evolution then occurs near a constant solution --- and not near a periodic wave. In turn once the slow part of the evolution is known one expects to be able to recover an asymptotic description of the full evolution including oscillations by using solutions of one-cell problems provided by profiles of periodic waves.

To implement this reduction --- even arguing on formal grounds --- one needs a way to derive averaged equations on one hand and on the other hand a way to translate initial data for the original system into effective initial data for the averaged system. The former may be achieved by inserting a slow/oscillatory \emph{ansatz}, that is, by studying a suitable family of solutions evolving from a family of well-prepared data. In contrast the main difficulty in answering the latter comes from the fact that one needs to understand what is the effective initial datum emerging from an ill-prepared datum. This problem may nevertheless be solved explicitly here by carrying out the analysis.

In this way we prove in \cite{JNRZ-conservation} that in the large-time regime the dynamics of parabolic systems near stable periodic waves follows at main order the slow modulation scenario that one may derive formally by some second-order version of the averaging method of Whitham. As a by-product this refined asymptotic description proves that the notion of space-modulated stability is sharp for generic parabolic systems and identifies what are the null conditions --- called phase uncoupling in \cite{JNRZ-conservation} --- that a given wave must satisfy to recover classical orbital stability from space-modulated stability. In \cite{R_linKdV} we also show that a similar scenario occurs for the linearized dynamics of \eqref{KdV}.

\section{Specific results}\label{s:results}

We now state precise results concerning \eqref{conservation} and \eqref{KdV}, as proved in \cite{JNRZ-conservation} and \cite{R_linKdV}.

\subsection{Periodic waves of parabolic systems}\label{s:parabolic}

We first give a precise definition of the notion of \emph{diffusive spectral stability}. To do so we fix a periodic wave of \eqref{conservation}, $\bU(t,x)=\ubU(\uk x+\uom t)$. Then we write \eqref{conservation} in the corresponding moving frame
\be\label{conservation-move}
\bU_t\,+\,\uom\bU_x\,+\,\uk\,(\bff(\bU))_x\,=\,\uk^2\,\bD\,\bU_{xx}
\ee
and denote by $L$ the generator of the corresponding linearized dynamics
$$
L\bV\ =\ -\,\uom\bV_x\,-\,\uk\,(\dD\bff(\ubU)(\bV))_x\,+\,\uk^2\,\bD\,\bV_{xx}\,.
$$
We shall say that the wave $\ubU$ is diffusively spectrally stable if the four conditions \eqref{D1}, \eqref{D2}, \eqref{D3} and \eqref{H} hold where

\be
\label{D1}
\begin{array}{r}
\sigma_{L^2(\R)}(L)\ \subset\ \left\{\ \lambda\in\C\ \middle|\ \textrm{Re}(\lambda)<0\ \right\}\ \cup\ \{0\}
\end{array}\tag{D1}
\ee

\be
\label{D2}
\begin{array}{r}
\textrm{There exists $\theta>0$ such that for any $\xi\in[-\pi,\pi)$\hspace{10em}}\\[1em] 
\sigma_{L^2_{per}((0,1))}(L_\xi)\ \subset\ \left\{\ \lambda\in\C\ \middle|\ \textrm{Re}(\lambda)\geq -\theta\,\xi^2\ \right\}\,.
\end{array}\tag{D2}
\ee

\be
\label{D3}
\begin{array}{l}
\textrm{$0$ is an eigenvalue of $L_0$ of algebraic dimension $d+1$.}
\end{array}\tag{D3}
\ee

\be
\label{H}
\begin{array}{l}
\textrm{First-order derivatives of critical curves}\\ 
\textrm{with respect to the Floquet exponent are distinct at $\xi=0$.}
\end{array}\tag{H}
\ee

It is not immediately clear that assumption \eqref{H} makes sense since $0$ is in general not semi-simple and generally speaking spectral perturbation of a Jordan block does not yield differentiable spectral curves. Nevertheless under assumption \eqref{D3} one may indeed prove that those spectral curves are differentiable at $\xi=0$ so that assumption \eqref{H} makes perfect sense\footnote{The most instructive proof of this leads to a connection at the spectral level with averaged dynamics. Yet alternatively classical arguments proving Rellich's theorem do show that this is also a consequence of \eqref{D1}.}. Moreover assumption \eqref{H} then implies that those critical curves are analytic with respect to the Floquet exponent. In particular, under assumptions \eqref{D1}, \eqref{D3} and \eqref{H}, assumption \eqref{D2} amounts to the non vanishing at $\xi=0$ of the second-order derivative of the real part of critical spectral curves with respect to the Floquet exponent. See \cite{Noble-Rodrigues} for detailed proofs and further discussions.

As emphasized by the choice in notation, assumption \eqref{H} plays a role distinct in the spectral stability and is not expected to be sharp for the obtained results. Indeed it plays a role similar to the assumption of strict hyperbolicity in the analysis of hyperbolic systems and could likely be replaced with something analogous to an assumption of symmetrizability. See \cite{R} for further discussions in this direction.

As is implicit in assumption \eqref{D3} one expects the family of periodic waves to be of dimension $d+2$ --- counting phase translations. Indeed when analyzing the profile system arising from \eqref{conservation}, it is immediate that all equations are conservative and hence may be integrated once by introducing $d$ constants of integration, besides the already-introduced phase speed. Once those $d+1$ parameters are fixed one does expect that around a given periodic orbit in the phase portrait of the reduced profile system no other periodic orbit co-exist. Hence in this case follows a local parametrization of periodic wave profiles by $d+2$ parameters. When a similar parametrization is available then \eqref{D3} really is an assumption of minimal multiplicity. In turn assuming \eqref{D3} one may prove that a parametrization by $d+2$ parameters holds and even identify wavenumber, averages over a period (of components corresponding to conservation laws) and phase shift as a set of suitable parameters. This essentially follows from an Evans function computation in the spirit of \cite{Serre}. See detailed discussions in \cite{JNRZ-conservation,Benzoni-Noble-Rodrigues,R}.

\begin{theorem}[\cite{JNRZ-conservation}]\label{th:stab-conservation}
Let $K\geq3$. For any $\ubU$ a periodic traveling wave profile of \eqref{conservation} that is diffusively spectrally stable there exist positive $\epsilon_0$ and $C$ such that for any $\bU_0$ such that 
$$
E_0:=\delta_{L^1(\R)\cap H^K(\R)}(\bU_0,\ubU)\ <\ \epsilon_0
$$ 
there exists $\bU$ a global classical solution of \eqref{conservation-move} with initial data $\bU_0$ such that for any $t\geq0$
$$
\delta_{H^K(\R)}(\bU(t),\ubU)\ \leq\ C\,E_0
%\delta_{L^1(\R)\cap H^3(\R)}(\bU_0,\ubU)
$$
and any $2\leq p\leq \infty$
$$
\delta_{L^p(\R)}(\bU(t),\ubU)\ \leq\ C\ (1+t)^{-\tfrac12\left(1-\tfrac1p\right)}\,E_0\,.
%\,\delta_{L^1(\R)\cap H^3(\R)}(\bU_0,\ubU)\,.
$$
\end{theorem}

The proof of the foregoing theorem also includes decay in $\delta_{H^K}$ but with non sharp decay rates. Indeed our proof is not optimized in regularity --- neither in decay of the high derivatives nor in the threshold $K\geq3$ --- but in robustness. In particular the proof does not rely directly on analytic semigroup techniques and applies almost word-by-word to quasilinear cases satisfying pointwise Kawashima conditions. See detailed discussion in \cite{JNRZ-conservation,R}. By similar techniques one does not expect to recover sharp decay rates in $\dot{W}^{k,p}$ for $\tfrac12K\leq k\leq K$ ; see \cite{Rodrigues-compressible} for instance.

The restriction $2\leq p\leq\infty$ has a distinct technical origin. It follows from our will to prove most of the required linear estimates on the spectral side by relying on Hausdorff-Young estimates as this leads to relatively simpler proofs than those relying on finer multiplier theorems or pointwise bounds of Green functions. See however \cite{Jung-pointwise-RD,Jung-pointwise-conservation} for examples of latter estimates.

In any case decay in $L^1$ should not be expected when the main decay mechanism is --- as here --- due to spatial dispersion of the solution. In this case, the trade-off encoded by estimates is localization against (regularity and) decay. The situation is in some sense parallel to the case where the decay mechanism is mixing and the trade is regularity against decay.

More generally one should not expect to be able to prove decay rates without loss --- in localization or in regularity or... --- when optimal decay rates are not of exponential type. At the linearized level the fact that this is indeed impossible is essentially a consequence of the Datko-Pazy theorem \cite[Theorem~3.1.5 \& Corollary 3.1.6]{vanNeerven}.

We now turn to large-time asymptotics. First we point out that assumption \eqref{D3} implies the existence of a family of nearby periodic waves with profiles $\ubU^{(\bM,k)}(\,\cdot\,+\phi)$ and phase speed $c(\bM,k)$, where $k$ is the corresponding wavenumber, $\phi$ is a phase shift and $\bM\in\R^d$ are averaged-values over a period
$$
\bM\ =\ \int_0^1 \ubU^{(\bM,k)}\,.
$$
As a consequence corresponding time frequency is $\omega(\bM,k)=-k\,c(\bM,k)$. Note carefully that what was $\ubU$ in Theorem~\ref{th:stab-conservation} is now $\ubU^{(\ubM,\uk)}$ for some $\ubM$ and that $\uom=\omega(\ubM,\uk)$. Following the extended version of the original strategy of Whitham \cite[Chapter~14]{Whitham} introduced in \cite{Noble-Rodrigues} --- see also \cite[Appendix~B]{JNRZ-conservation} ---, one may derive arguing on formal grounds that two-scale slow/oscillating solutions of \eqref{conservation} should behave as
\be\label{slow-modulation}
(t,x)\ \mapsto\ \ubU^{(\cM,\kappa)(t,x)}\ \left(\Psi(t,x)\right)
\ee
with $\kappa=\d_x\Psi$, where $(\cM,\kappa)$ are slowly evolving according to some system
\ba\label{whitham}
\cM_t+(\bF(\cM,\kappa))_x &=(d_{11}(\cM,\kappa) \cM_x + d_{12}(\cM,\kappa)\kappa_x)_x\,,\\
\kappa_t-(\omega(\cM,\kappa))_x&=(d_{21}(\cM,\kappa)\cM_x+ d_{22}(\cM,\kappa)\kappa_x)_x\,.
\ea
Here $\bF$ is simply the averaged flux
$$
\bF(\bM,k)\ =\ \int_0^1\bff(\ubU^{(\bM,k)}(x))\,\dD x
$$
and $d_{i,j}$ are determined in a more complicated way as usually expected from higher-order corrections in averaging processes.

Though we do not expound the full derivation of \eqref{whitham} we emphasize now two key-points of this derivation. The first one is that the last equation of \eqref{whitham} is first obtained as
\be\label{psieq}
\Psi_t= \omega(\cM,\kappa )+ d_{21}(\cM,\kappa)\cM_x + d_{22}(\cM,\kappa)\kappa_x
\ee
that is then differentiated with respect to the space variable. It is important to note that a completely similar scenario occurs at the spectral level, in particular when desingularizing the Jordan block structure, as the generalized kernel is spanned by $\ubU_x$, $\d_{\bM_1}\ubU$, $\cdots$, $\d_{\bM_d}\ubU$, but one needs to unravel the role of $\d_k\ubU$ by transforming a phase-like component in a wavenumber-like component through a suitable multiplication by $\iD\xi$ ; see \cite{Noble-Rodrigues,Benzoni-Noble-Rodrigues,KR}. The second point is that for our purpose of describing to main order the large-time dynamics near $(\ubM,\uk)$ system~\eqref{whitham} is far from being uniquely determined but there is a canonical choice that also ensures that though the system has been derived by fitting slow evolutions --- a low-Fourier requirement --- the obtained system is a well-posed hyperbolic-parabolic system --- a high-Fourier feature. See \cite{Noble-Rodrigues} and \cite[Appendix~B]{JNRZ-conservation}.

To compare directly to solutions of \eqref{conservation-move} (and not \eqref{conservation}) we now denote by $(W)$ and $(W)_{phase}$ system \eqref{whitham} and equation \eqref{psieq} written in the frame of $\ubU^{(\ubM,\uk)}$.

\begin{theorem}[\cite{JNRZ-conservation}]\label{th:asymptotic-conservation}
Let $K\geq4$. For any $\ubU$ a periodic traveling wave profile of \eqref{conservation} that is diffusively spectrally stable there exist positive $\epsilon_0$ and $C$ such that for any $\bU_0$ such that there exists $\psi_0$ such that
$$
E_0:=\|\bU_0\circ(\Id-\psi_0)-\ubU\|_{L^1(\R)\cap H^K(\R)}\ +\ \|\d_x\psi_0\|_{L^1(\R)\cap H^K(\R)}\ <\ \epsilon_0
$$ 
there exist $\bU$ a global classical solution of \eqref{conservation-move} with initial data $\bU_0$, a local phase shift $\psi$ with initial data $\psi_0$ and local averages $\bM$ such that for any $t\geq0$ and any $2\leq p\leq \infty$
\be\label{est:slow-modulation}
\begin{array}{rl}
\|\bU(t,\cdot-\psi(t,\cdot))\ -\ \ubU^{\left(\ubM+\bM(t,\cdot),\frac{\uk}{(1-\psi_x(t,\cdot))}\right)}(\,\cdot\,)\|_{L^p(\R)}%\\[1em]
&\leq\ C\,E_0\ \ln(2+t)\ (1+t)^{-\frac{3}{4}}\\[1em]
\|(\bM,\uk\,\psi_x)(t,\cdot) \|_{L^p(\R)}
&\leq\ C\,E_0\ (1+t)^{-\frac{1}{2}(1-1/p)}\\[1em]
\|\psi(t,\cdot)-\tfrac{\psi_0(-\infty)+\psi_0(+\infty)}{2}\|_{L^\infty(\R)}&\leq\ C\,E_0\,;
\end{array}
\ee
and setting $\Psi(t,\cdot)=(\Id-\psi(t,\cdot))^{-1}$, $\kappa=\uk\,\d_x\Psi$, $\cM(t,\cdot)=(\ubM+\bM(t,\cdot))\circ\Psi(t,\cdot)$, and letting $(\cM_W,\kappa_W)$ solve $(W)$ with initial data 
\be\label{data}
\begin{array}{rcl}
\kappa_W(0,\cdot)&=&\uk\,\d_x\Psi(0,\cdot)\,,\\[0.5em]
\cM_W(0,\cdot)&=&\ubM+ \bU_0-\ubU\circ\Psi(0,\cdot)%\\
%&+&
+\left(\dfrac{1}{\d_x\Psi(0,\cdot)}-1\right)\left(\ubU\circ\Psi(0,\cdot)-\ubM\right)\,,
\end{array}
\ee
we have
\be
\|(\cM,\kappa)(t,\cdot)-(\cM_W,\kappa_W)(t,\cdot) \|_{L^p(\R)}
\ \leq\ C\,E_0\ (1+t)^{-\frac{1}{2}(1-1/p)-\frac{1}{2}+\eta}\,;%\\[0.5em]
\ee
at last recovering $\Psi_W$ from $(\cM_W,\kappa_W)$ though $(W)_{phase}$ with $\Psi_W(0,\cdot)\ =\ \Psi(0,\cdot)$
\be
\|\Psi(t,\cdot) -\Psi_W(t,\cdot) \|_{L^p(\R)}
\ \leq\ C\,E_0\ (1+t)^{-\frac{1}{2}(1-1/p)+\eta}.%\\[0.5em]
\ee
\end{theorem}

The upshot of estimates \eqref{est:slow-modulation} is that $\bU$ is indeed very well approximated by modulating $\ubU$ in all its parameters, showing that asymptotically $\bU$ essentially takes the form \eqref{slow-modulation}. The remaining estimates prove that the evolution of slow parameters is indeed well-captured by \eqref{whitham}. Observe that in order to achieve relevant comparisons one needs to stick to the description \eqref{slow-modulation} and effectively perform all modulations on $\ubU$ hence to invert $\Id-\psi(t,\cdot)$. 

Observe also how \eqref{data} prescribes data in a non trivial way. In particular the last term in $\cM_W(0,\cdot)$ incorporates contributions of the high frequencies of $\psi_0$ and hence cannot be included in a formal slow/oscillatory \emph{ansatz}.

Relying then on a fine description of near constant dynamics --- see \cite{Liu_Zeng} and \cite[Appendices~B \&~C]{JNRZ-conservation} ---, one may then use Theorem~\ref{th:asymptotic-conservation} to obtain a very precise description of the dynamics in terms of weakly interacting diffusion-waves in local parameters. This yields that decay rates in Theorem~\ref{th:stab-conservation} are sharp for general systems but also enables us to identify various sharp cancellation conditions, coined as phase uncoupling conditions in \cite{JNRZ-conservation}, that leads when $\psi_0$ is constant (or when $\d_x\psi_0$ is mean-free) to extra decay hence to classical (and not space-modulated) asymptotic orbital stability and to an asymptotically linear behavior.

\subsection{Cnoidal waves of the Korteweg--de Vries equation}\label{s:KdV}

We now discuss similar results for the \emph{linearized} dynamics of \eqref{KdV} near cnoidal waves. 

First select such a cnoidal wave $\uU$ and consider the operator $L$ introduced in \eqref{linop}. One readily checks that $L$ does generates a group $(S(t))_{t\in\R}$ on $L^2(\R)$ (or on $H^s(\R)$ with $s\in\N$, or...) by using energy estimates which also lead to the crude bound
$$
\|S(t)\|_{L^2(\R)\to L^2(\R)}\ \lesssim\ e^{|t|\,\tfrac{\uk}{2}\|(\uU_x)_-\|_{L^\infty(\R)}}\,,\quad t\in\R\,.
$$
That energy estimates yield such a poor bound is a manifestation of the fact that $L$ is not a normal operator.

Nevertheless space-modulated bounded linear stability holds in $H^s(\R)$ for any $s\in\N$.

\begin{theorem}\label{th:stab-KdV}
For any cnoidal wave of \eqref{KdV} and any $s\in\N$, there exists $C$ such that for any $W_0$ and any time $t\in\R$
$$
N_{H^s(\R)}(S(t)\,W_0)\ \leq\ C\, N_{H^s(\R)}(W_0)
$$
where $(S(t))_{t\in\R}$ is the group of operator solutions of the linearized dynamics.
\end{theorem}

The foregoing theorem is non trivial even when $s=0$. Yet there does not seem to be an easy way to deduce the general case from the case $s=0$. Since the emphasis here is in developing robust techniques rather than in the result itself, we now give some hints about its proof\footnote{In contrast we have chosen not to comment proofs of Theorems~\ref{th:stab-conservation} \&~\ref{th:asymptotic-conservation} since they are already sketched in the introduction of \cite{JNRZ-conservation} itself, further commented in \cite{R} and one may gain some insight about their key-features by looking first at the analysis of a simpler sub-case treated in \cite{JNRZ-RD1,JNRZ-RD2}. Incidentally we mention that this sub-case also exhibits the very specific feature that at the averaged level it is asymptotically self-similar and hence is amenable to a resolution by the renormalization techniques initially used by Schneider; see \cite{SSSU}.}.

The starting point is that one knows that the spectrum of $L$ lies on the imaginary axis --- see for instance \cite{Bottman-Deconinck,R_linKdV}. Yet by itself this is insufficient to derive any suitable bound. To go beyond one may build for each Floquet $\xi$ a spectral decomposition of $L_\xi$. This is essentially a consequence of the theory pioneered by Keldy\v{s} \cite{Keldysh-1,Keldysh-2} --- see \cite{Markus,Yakubov,Gohberg-Goldberg-Kaashoek-1} --- on spectral decompositions of relatively compact perturbations of self-adjoint operators with compact resolvents. It relies on the fact that as far as spectrum at infinity is concerned $L_\xi$ may be thought as a perturbation of the skew-adjoint $-\uk^3(\d_x+\iD\xi)^3$. With this decomposition in hands and knowing that except zero as an eigenvalue of $L_0$ all Floquet eigenvalues are simple, after a desingularization of the zero eigenvalue in a way similar to what was required in the analysis of parabolic cases --- see \cite{Benzoni-Noble-Rodrigues,R_linKdV} --- one may then define Sobolev-like norms $\|\cdot\|_{X^s(\R)}$ for which the evolution generated by $L$ is unitary. What is left then is to show that $\|\cdot\|_{X^s(\R)}$ is indeed equivalent to $N_{H^s(\R)}$. This is actually the core of the proof but relies, as above arguments, mostly on fine analyses of large and small spectrum.

While some proofs are easier to carry out by using explicit formulas derived from Lax pair representation, most of them may be obtained in a robust way. Notable exceptions are the knowledge that the spectrum indeed lies on the imaginary axis --- a minimal requirement for linear stability --- and that non zero Floquet eigenvalues are indeed simple.

In contrast in order to derive asymptotic stability we do use a huge wealth of information. Indeed our argument uses a non degeneracy assumption on the Floquet parametrization of the spectrum similar to \eqref{D2} but here the whole spectrum is critical hence instead of requiring such a condition at one point, zero, as in \eqref{D2} we need such a piece of information everywhere. To be a bit less pessimistic let us precise that actually both at infinity and near zero the required pieces of information may be obtained by robust techniques so that what remain is really a sign control on smooth functions defined on compact intervals\footnote{This is typically the kind of verification amenable to computer-assisted proofs relying on interval arithmetics as in \cite{Barker}.}.

\begin{figure}
\begin{center} 
\includegraphics[scale=0.5]{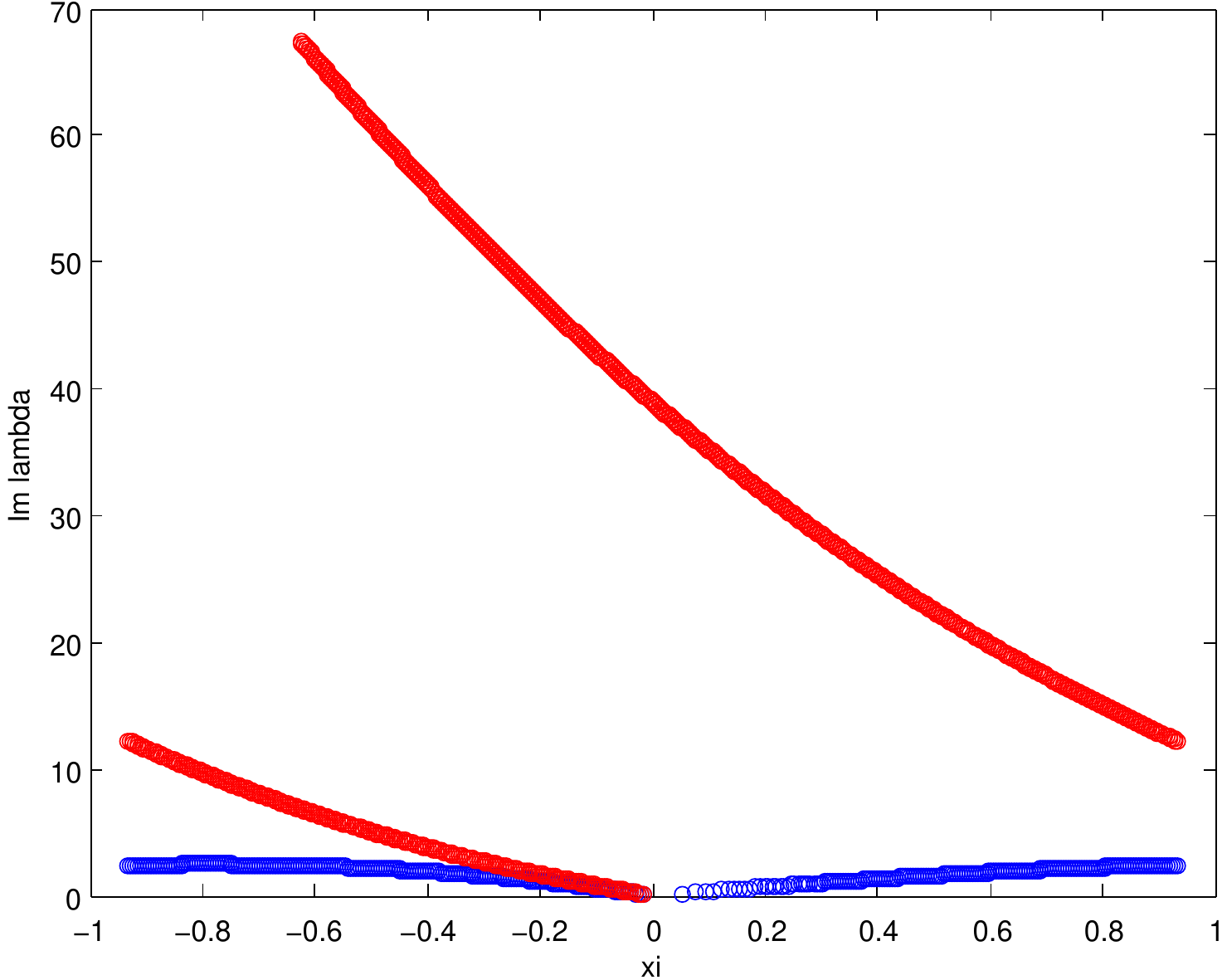}
\end{center}
\caption{For a typical cnoidal wave, part of the spectrum of $L$ with positive imaginary part : $\textrm{Im}(\lambda)$ \emph{vs} $\xi$ such that $\lambda\in\sigma_{L^2_{per}((0,1))}(L_\xi)$. There are infinitely many red branches going to infinity.}
\label{spectrum-sketch}
\end{figure}

Before giving a precise statement of this non degeneracy condition we give more details concerning the structure of the spectrum of $L$. First, real and Hamiltonian symmetries yield that $\lambda\in\sigma_{L^2_{per}((0,1))}(L_\xi)$ if and only if $-\lambda\in\sigma_{L^2_{per}((0,1))}(L_{-\xi})$. As a by-product of the proof of Theorem~\ref{th:stab-KdV} we also know on one side that the large spectrum consists in two symmetric branches that by gluing together various Floquet parametrizations may be parametrized by two symmetric semi-unbounded intervals, with a parametrization equivalent to $\xi\mapsto -(\iD\uk\xi)^3$ and whose third order derivative of the imaginary part converges to $6\uk$ and on the other side that near zero there are three branches of spectrum going through zero with second-order derivative vanishing at zero but third order derivatives non zero. By the techniques used in \cite{Bottman-Deconinck} one may actually derive a complete parametrization of $(\lambda,\xi)$ in terms of an auxiliary parameter --- a spectral Lax parameter --- and conclude that this picture extends as follows: the spectrum may be split in two parts, one infinite line covering the imaginary axis, and a symmetric loop covering twice a neighborhood of zero in the imaginary axis\footnote{The presence of this triply-covered area near zero does not contradict the simplicity of Floquet eigenvalues as the passage at the same point of the spectrum occurs at distinct Floquet parameters.}; see Figure~\ref{spectrum-sketch}. The picture is of course consistent with descriptions obtained by more robust techniques in asymptotic regimes. In particular, in the solitary wave limit \cite{Gardner-large-period} the line coincides with the essential spectrum of solitary waves and the loop arises from the embedded eigenvalue zero.

\medskip

Despite the relatively explicit description of the spectrum of $L$, up to now the author has not been able to prove the following condition
\be
\label{A}
\begin{array}{l}
\textrm{Along the line the third-order derivative with respect to Floquet exponent}\\
\textrm{does not vanish. Along the loop the second-order derivative}\\ 
\textrm{with respect to Floquet exponent does not vanish except at zero.}
\end{array}\tag{A}
\ee
However for all cnoidal waves picked randomly by the author an eye inspection of corresponding graphs seems to show that condition~\eqref{A} is satisfied.

\begin{theorem}\label{th:asymptotic-KdV}
For any cnoidal wave of \eqref{KdV} satisfying \eqref{A}, there exists $C$ such that for any $W_0$ and any time $t\in\R^*$
$$
N_{L^\infty(\R)}(S(t)\,W_0)\ \leq\ C\,|t|^{-1/3}\,N_{L^1(\R)}(W_0)
$$
and
$$
N_{L^\infty(\R)}(S(t)\,W_0)\ \leq\ C\,(1+|t|)^{-1/3}\,N_{L^1(\R)\cap H^1(\R)}(W_0)
$$
where $(S(t))_{t\in\R}$ is the group of operator solutions of the linearized dynamics.
\end{theorem}

In contrast with most of results introduced so far our proof of Theorem~\ref{th:asymptotic-KdV} is not carried out on the Bloch side but proceeds by pointwise estimates of suitable Green functions through oscillatory integral estimates. As in the parabolic case one may actually go further and prove at this linearized level that the large-time dynamics is at main order of slow-modulation type and that the involved local parameters evolve essentially according to a third-order hyperbolic-dispersive system. For the sake of brevety we do not detail this here but add two comments. First, as the reader may have already deduced from the structure of the spectrum of $L$, the dimension of the family of cnoidal waves is $4$ (counting phase shift), the extra parameter --- compared to phase shift, wavenumber, average of $U$ ---, which arises from the Hamiltonian structure and invariance by space translation, may be chosen as the average of the Benjamin impulse; see \cite{Benzoni-Noble-Rodrigues,Benzoni-Mietka-Rodrigues}. Second, we stress that to capture the large-time effects encoded by the form of spectral curves near zero one does need to include a third-order correction to the classical first-order Whitham system.

At last we warn the reader that it is very unlikely that decay rates expounded here could serve directly in a nonlinear large-time study since they seem by far too slow for such a purpose, and that the lack of normality of $L$ prevents a direct use of the classical $T-T^*$ argument to deduce from dispersive estimates of Theorem~\ref{th:asymptotic-KdV} bounds of Strichartz type.
 
%%%%%%%%%%%%%%%%%%%%%%%%%%%%%%%%%%%%%%%%%%%%%%%%%%%%%% 
%        bibliography                                %
%%%%%%%%%%%%%%%%%%%%%%%%%%%%%%%%%%%%%%%%%%%%%%%%%%%%%%

\bibliographystyle{abbrv}
\bibliography{Ref} 

\end{document}